\documentclass{IEEEtran4PSCC}
\usepackage{hyperref}
\usepackage{cite}
\ifCLASSINFOpdf
   \usepackage[pdftex]{graphicx}
\else
   \usepackage[dvips]{graphicx}
\fi
\usepackage[cmex10]{amsmath}
\usepackage{setspace}
\usepackage{url}
\usepackage{xcolor}
\usepackage{booktabs}
\newcommand{\mc}{\mathcal}
\newcommand{\ts}{\text}
\makeatletter
\let\old@ps@headings\ps@headings
\let\old@ps@IEEEtitlepagestyle\ps@IEEEtitlepagestyle
\def\psccfooter#1{%
    \def\ps@headings{%
        \old@ps@headings%
        \def\@oddfoot{\strut\hfill#1\hfill\strut}%
        \def\@evenfoot{\strut\hfill#1\hfill\strut}%
    }%
    \def\ps@IEEEtitlepagestyle{%
        \old@ps@IEEEtitlepagestyle%
        \def\@oddfoot{\strut\hfill#1\hfill\strut}%
        \def\@evenfoot{\strut\hfill#1\hfill\strut}%
    }%
    \ps@headings%
}
\makeatother

\psccfooter{%
  \parbox{\textwidth}{\hrulefill \\ \centering \thepage}%
}

\begin{document}
\title{Modelling Renewable Curtailment and Constraints in Ireland's Electricity System}

\author{
\IEEEauthorblockN{Richard Nayer, Sam Hodges, Waqquas Bukhsh}
\IEEEauthorblockA{University of Strathclyde\\Glasgow, Scotland (UK)\\
\{richard.nayer,sam.hodges,waqquas.bukhsh\}@strath.ac.uk}
\and
\IEEEauthorblockN{Claver Chitambo, Chanura Wijeratne}
\IEEEauthorblockA{Renewable Energy Solutions \\
Kings Langley, UK\\
\{claver.chitambo,chanura.wijeratne\}@res-group.com}
}

\maketitle

\begin{abstract}
This paper describes the electricity markets and operational processes in the Irish power system and translates them into a Mixed Integer Linear Programming (MILP) model. The model is designed to estimate renewable generation \emph{Curtailment} and \emph{Constraint}. A full mathematical formulation is presented and tested on both a simple example and a realistic model of the Irish transmission network. The results show that the proposed formulation works as intended, capturing curtailment and constraint effects, while also highlighting model limitations and directions for future improvement.
\end{abstract}


% % Use this to place sponsorships
% \thanksto{\noindent Submitted to the 24th Power Systems Computation Conference (PSCC 2026).}

\section{Introduction}
In market based power systems, generators have traditionally received \emph{firm} connection agreements from Distribution and Transmission Network Operators (DNOs/TNOs), ensuring uninterrupted network access. However, the rapid build out of renewable generators at both transmission and distribution level has outpaced network reinforcement, leading to long lead-times and high costs for securing firm connections. To address this, \emph{non-firm} connection agreements have become increasingly common, under which grid access may be curtailed without compensation when network constraints or system security requirements arise. Such arrangements expose renewable generators to revenue uncertainty. A review of international practices in \cite{flex_connections_review} highlights the growing importance of forecasting export reductions, which is critical for developers' investment decisions.

\begin{figure}
    \centering
    \includegraphics[width=1\linewidth]{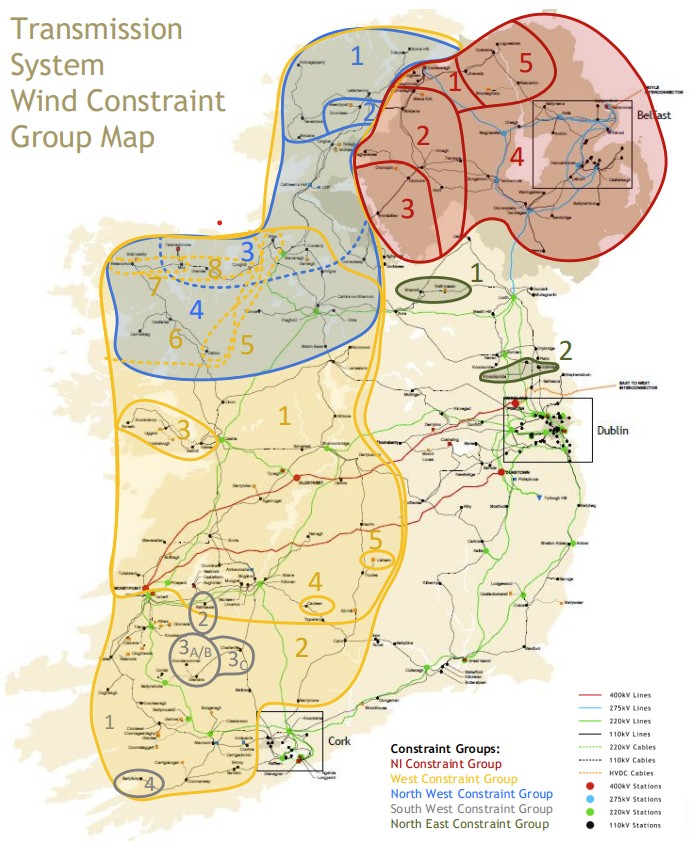}
    \caption{Ireland \emph{Constraint} Groups (While the figure states "Wind Constraint Group Map", solar generators will also be constrained within these pro-rata groups) \cite{EIRGRID_A}.}
    \label{fig:ConstraintGroups}
\end{figure}

The process by which non-firm connections and export limits are managed differs across power systems and operators~\cite{flex_connections_review}. For example, in Great Britain, DNOs apply a Last-In-First-Out (LIFO) principle, where the most recent connections are curtailed first in real time~\cite{flex_connections_review}. On the Island of Ireland, the TNOs use a pro-rata approach, where curtailments are shared proportionally among renewable generators~\cite{flex_connections_review}. In South Australia, the DNOs apply Dynamic Operating Envelopes (DOE), with capacity limits updated ahead of each dispatch interval~\cite{flex_connections_review}. This paper focuses on the island of Ireland, where pro-rata reductions can be applied across all non-synchronous generators on the transmission network to address system security requirements, or to specific \emph{constraint groups} (which may overlap) to address local transmission constraints. Figure~\ref{fig:ConstraintGroups} presents a schematic of the high-voltage transmission network in Ireland, with the main constraint groups shown geographically~\cite{EIRGRID_A}. These groups are not static; they change over time depending on the expected operating conditions, which in turn are influenced by the level of non-synchronous generation, the network state, cross-border imports and exports, and the availability of wind and solar resources.

%<<<<<<<<<<<<<<<<<<<
%SECTION CHANGE
%>>>>>>>>>>>>>>>>>>>

\section{Literature Review \& Contributions}\label{sec:lit_review}
The use of mathematical modelling for power system operation and planning is well established in both industry and academia, with a wealth of associated literature \cite{FERC_paper_1, KAYA_fiftyyearsofpowersystemsoptimization, motta_surveypoweropt}. Several open source power system optimisation tools have been developed \cite{oats, PYPSA, pandapower.2018}, many of which can implicitly model generator curtailment. However, most of these tools focus on generic constraints and do not capture system-specific curtailment rules, such as pro-rata curtailment in overlapping constraint groups. In this paper, we build on a modified version of OATS \cite{oats}, which has previously been applied to curtailment analysis in \cite{hawker_hydrogencurtailment}. 

While there is some literature related to pro-rata curtailment, there is limited research that reflects the detailed curtailment and constraint approaches used in Ireland. EirGrid, the Irish Transmission System Operator, regularly assess renewable curtailment using a commercial platform called PLEXOS~\cite{Eirgrid_balancingprinciples}, but the exact model formulation is not available. In \cite{6866255}, a MATPOWER model performs optimal power flow (OPF) assessments of LIFO, Rota and Pro-Rata curtailment policies. The method assumed a single pro-rata group, whereas in Ireland generators can belong to multiple, and potentially overlapping, pro-rata groups. A different approach in \cite{ANDONI2017174} used a game theoretic model with statistical distributions of generation and demand to compare curtailment under LIFO, pro-rata, and rota policies. However, the approach does not explicitly model the network and was limited to two demand regions. The recent work in \cite{10253224} analysed dispatch instruction data from 2020-2021 to estimate curtailment volumes in Ireland. The study found that most system level curtailment resulted from minimum unit run requirements (MUON), with the remainder driven by the System Non-Synchronous Penetration (SNSP) limit of 75\%. SNSP is defined by EirGrid/SONI as~\cite{EirGridSONI_SNSPCalc}:
\begin{center}
 \textit{``the ratio of the real-time MW contribution from non-synchronous generation (including BESS discharges and interconnector imports) to demand (including BESS charges and interconnector exports)"}.   
\end{center}
  
To best of authors' knowledge, no published work provides an Irish electricity system-specific definition of curtailment and constraints, a step-by-step representation of how these are applied in practice, or a mathematical model that captures the overlapping and dynamic nature of constraint groups. This paper addresses that gap. Our contributions are twofold:

\begin{itemize}
    \item Elicitation and formalisation of the high-level requirements of the curtailment and constraint framework in Irish electricity system.
    \item Development of a simple yet representative market and network model for Ireland that can estimate curtailment and constraint volumes while accounting for the complexities of overlapping constraint groups.
\end{itemize}

The remainder of the paper is organised as follows. Section~\ref{sec:IrelandSystemOperation} provides an overview of the Irish electricity market and introduces definitions of curtailment and constraint. Section~\ref{sec:model} presents the mathematical model. Section~\ref{sec:test} demonstrates the results using both an illustrative test system and a representative Irish system. The paper concludes in Section~\ref{sec:conclusion}.

%<<<<<<<<<<<<<<<<<<<
%SECTION CHANGE
%>>>>>>>>>>>>>>>>>>>
\section{Island of Ireland Electricity System Operation}\label{sec:IrelandSystemOperation}
In Ireland, the terms \emph{Curtailment} and \emph{Constraint} have specific definitions, which are emphasised here to distinguish them from more general uses. The electricity transmission network across the Republic of Ireland (ROI) and Northern Ireland (NI) is jointly operated by the two TSOs, Eirgrid and SONI. Together they also operate the Single Electricity Market Operator (SEMO) for the all-island market. As part of the EU internal market, they also act as Nominated Electricity Market Operators (NEMOs) through the SEMO Power Exchange (SEMOpx). SEMOpx operate the following three wholesale markets ~\cite{SEMO_MarketProcedures}.

\begin{itemize}
    \item The Day-Ahead Market (DAM), which sets day-ahead price in Ireland and interconnection flows determined by the EUPHEMIA algorithm~\cite{Euphemia_description} .
    \item The Intraday Auctions (IDA1, IDA2, IDA3), which allow adjustments closer to real time. IDA1 and IDA2 are coupled with Great Britain, while IDA3 is local to Ireland.
    \item The Intraday Continuous Market (IDC), which is also local to Ireland. 
\end{itemize}

\begin{figure}[t]
\includegraphics[width=1\linewidth]{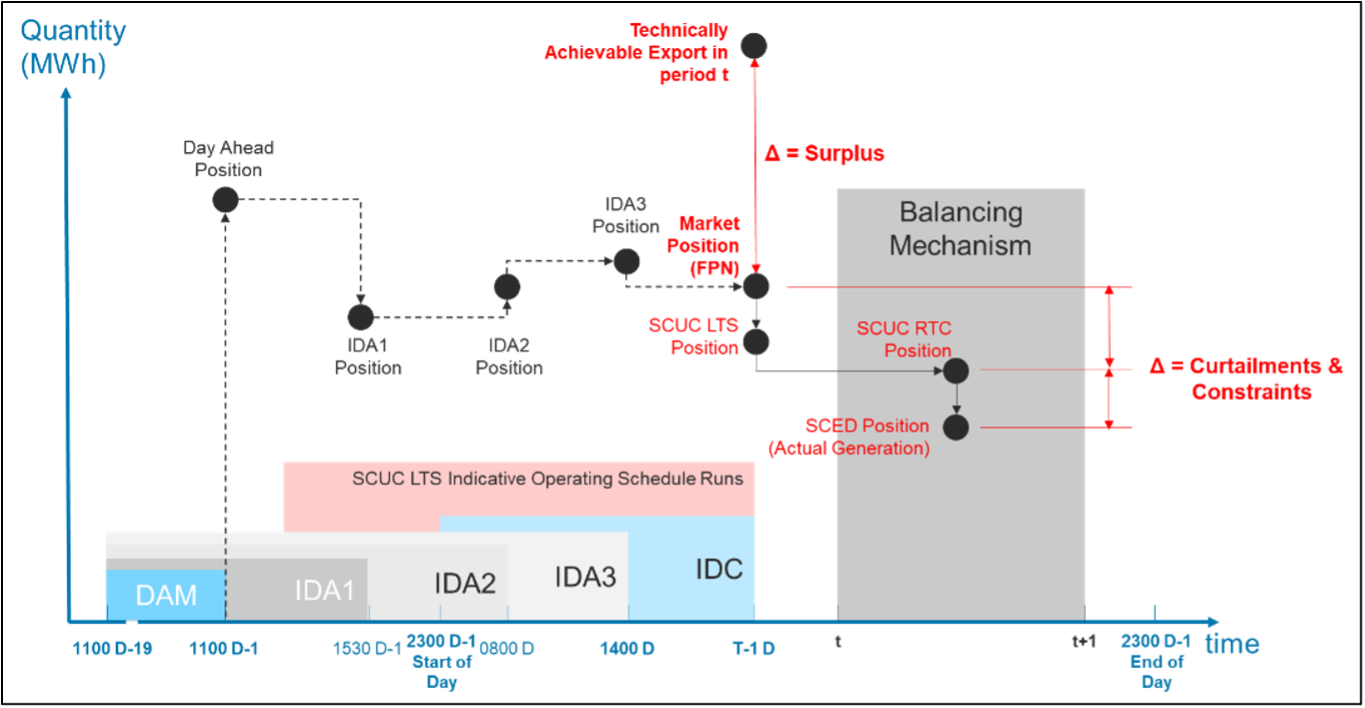}
\caption{All-Ireland Market Formation Sketch}
\label{fig:market}
\end{figure}

Imbalances between traded and actual generation are settled in the Balancing Mechanism, managed by SEMO. The DAM clears in 1-hour trading periods, while the intraday markets and balancing market operate in 30-minute periods (with 5-minute pricing in the balancing market). The last intraday market to close is the IDC, which stops trading one hour before delivery, known as Gate Closure 2 (GC2). At this point, generators submit their Final Physical Notification (FPN) of expected output. It should be noted that for the majority of renewable generators an FPN is not required to be submitted, and the forecast value will be taken. Figure~\ref{fig:market}) presents an example of how the position of a participant in the markets could vary with time.

After the gate closure, the TSOs run a Security Constrained Unit Commitment (SCUC) to produce indicative schedules for up to 48 hours ahead. This process, known as Long Term Scheduling (LTS), gives
market participants time to adjust their positions. Closer to real-time, further SCUC runs are performed in the Rea-Time Commitment (RTC), followed by Security Constrained Economic Dispatch (SCED), which informs dispatch instructions issued from the control room \cite{EirGridSONI_windprocess, EirGridSONI_dispatchprocess}. Deviations from FPNs are priced in the balancing market. The balancing market follows three main principles, in the following order of priority~\cite{Eirgrid_balancingprinciples}:

\begin{enumerate}
\item Ensuring operational security (by performing curtailment if required);
\item Maximise the use of priority-dispatch generation (mainly renewables); and,
\item Ensure the efficient operation of the SEM.
\end{enumerate}

\begin{figure}[t]
    \centering
    \includegraphics[width=1\linewidth]{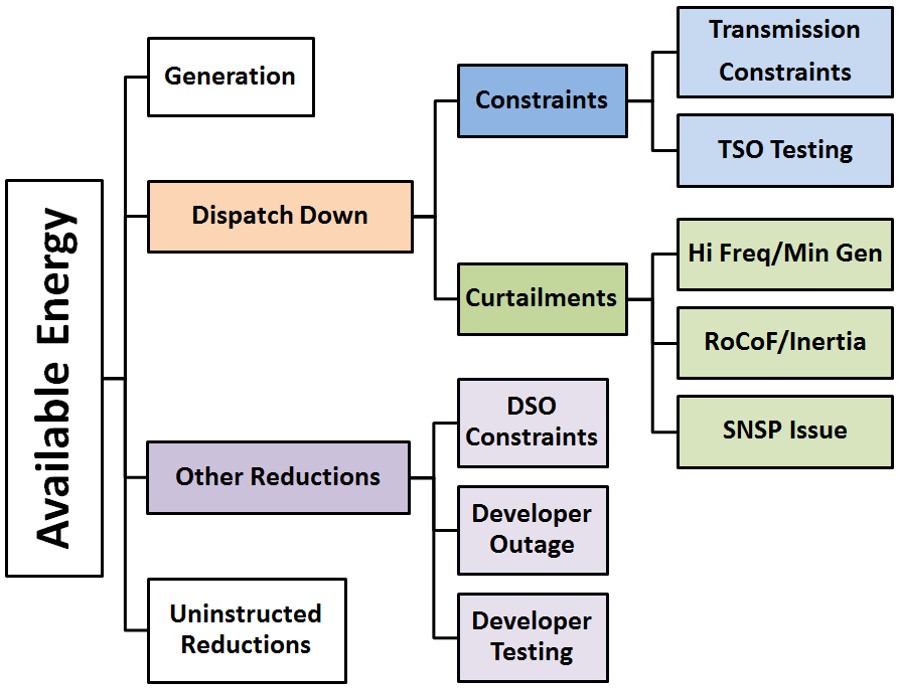}
    \caption{\emph{Curtailment} \& \emph{Constraint} Types \cite{EirGrid_2024_dd_report}}
    \label{fig:constraintcurtailment_types}
\end{figure}

The TSOs classify dispatch actions on non-synchronous generators (wind, solar, interconnectors) into \emph{Curtailment}, \emph{Constraint} or Other (see Fig.~\ref{fig:constraintcurtailment_types}, as described in the following sections. These volumes are reported publicly by SEMO. An additional measure of interest is Surplus, which reflects oversupply of renewables relative to demand~\cite{Eirgrid_ECP2.4_Methodology}. 

\subsection{Surplus}
Surplus occurs when the potential to generate power from renewables in a period exceeds the power demand on the system. It can be considered to be the difference between the technically achievable export from a generator (e.g. for renewables the forecast generation based on the weather), relative to their position achieved in the market (noting renewables are not required to submit FPNs). 

\subsection{Curtailment}
\emph{Curtailment} is dispatch down from the FPN to ensure system security. The reasons for curtailment are weekly published~\cite{SEMO_weeklyconstraintupdate}. Common reasons include the SNSP limit, minimum inertia or frequency requirements, and the need to keep a minimum number of conventional units online.

\subsection{Constraint}
\emph{Constraint} is the dispatch down (relative to the FPN) by the control room of renewable generators due to local transmission constraints. Renewables are all allocated to \emph{Constraint} groups (as shown in Fig.\ref{fig:ConstraintGroups}), which are defined by modelling to identify which generators contributed to particular network constraints~\cite{EIRGRID_A}.

Renewable generation is currently prioritised to avoid dispatch down, as required under EU and Irish regulations~\cite{Eirgrid_balancingprinciples}. In SCUC and SCED runs, their bid prices are replaced with negative values to ensure they are scheduled ahead of conventional plants. A 2019 decision \cite{SEMO_SEM-21-027} stated that new renewable generators would lose priority dispatch status, but in practice this has not yet been fully implemented~\cite{SEMO_SEM-22-009}.
%<<<<<<<<<<<<<<<<<<<
%SECTION CHANGE
%>>>>>>>>>>>>>>>>>>>

\section{Model Description}\label{sec:model}
\begin{figure}[h]
    \centering
    \includegraphics[width=1\linewidth]{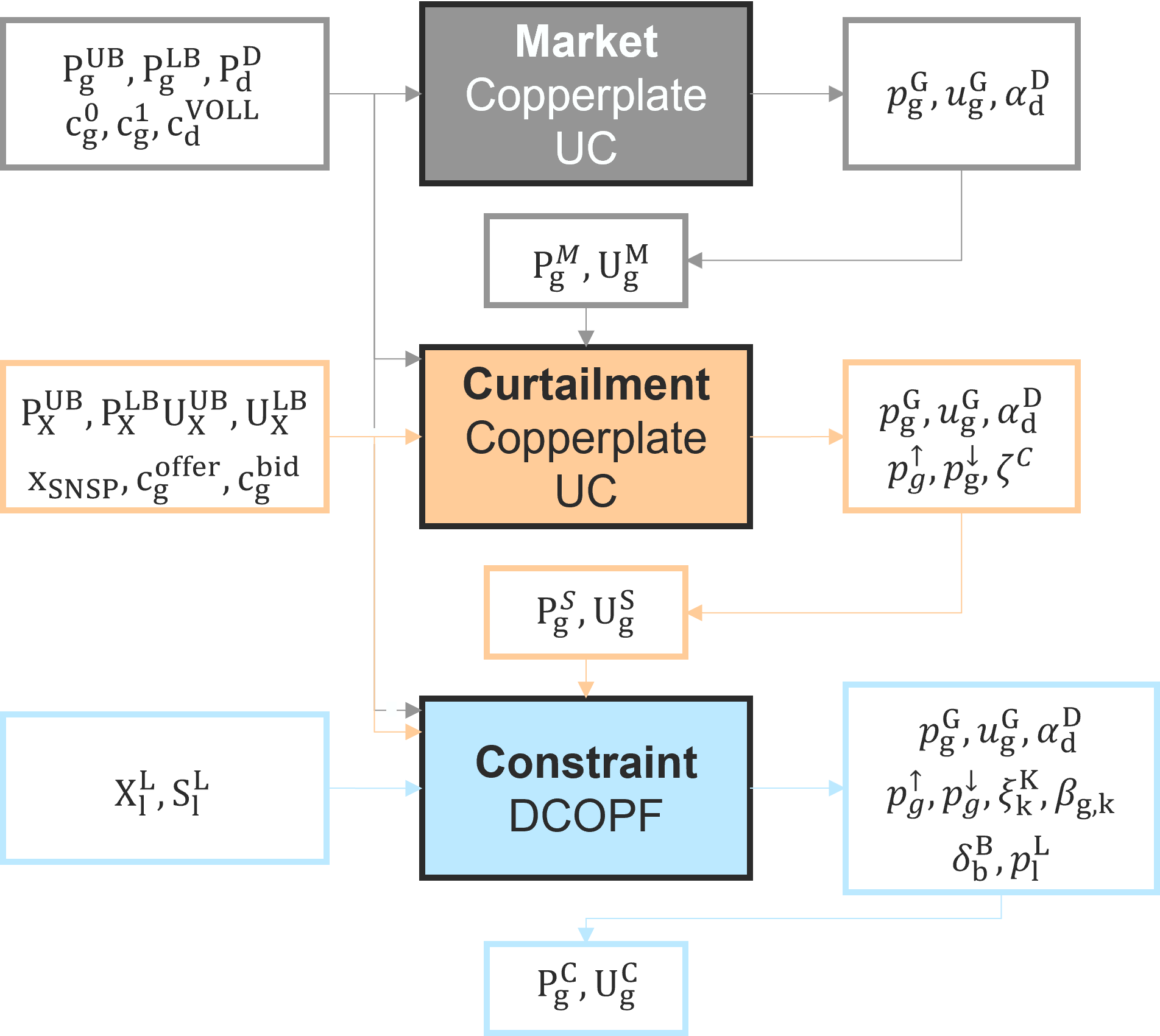}
    \caption{Model Flow Diagram (each snapshot). Parameters are shown on the left, with variables on the right.}
    \label{fig:model_flow}
\end{figure}

The model significantly reworks and extends upon an existing modelling framework developed at the University of Strathclyde, OATS \cite{oats}. The high level requirements and approach of the model (as summarised in Fig.~\ref{fig:model_flow}) are as follows. Note that the model is not time coupled, and ran as a series of snapshots:
\begin{enumerate}
    \item Market: Determine the ex-ante market position for each generator. Copper-Plate Unit Commitment (UC) is ran as surrogate to provide market derived generation set-points.
    \item Curtailed: Determine the volume of SNSP and MUON curtailment for renewable generation. Additional constraints are added to the Copper-Plate UC to incorporate market outcomes, and replicate the constraints applied by EirGrid to secure the system (e.g. SNSP, MUON).
    \item Constrained: Determine the volume of renewable constrain dispatch down. The model is reformulated as a UC DCOPF to account for local transmission constraints, and determined further redispatch.
\end{enumerate}

\subsection*{Notation}
\setstretch{1.2}
\subsubsection*{Sets}

\begin{description}[\IEEEsetlabelwidth{$\ts{P}^{\ts{LB}}_g,\ts{P}^{\ts{UB}}_g$}\IEEEusemathlabelsep]
\item[$\mc{B}$] Buses, indexed by $b$
\item[$\mc{G}$] Generators, indexed by $g$
\item[$\mc{D}$] Demands, indexed by $d$
\item[$\mc{L}$] Lines, indexed by $l=(b,b')$
\item[$\mc{C}$] \emph{Constrained} groups, indexed by $c$
\item[$\mc{K}_i$] $i^{th}$ \emph{Constraint} group, indexed by $k$
\end{description}

\subsubsection*{Parameters}
\begin{description}[\IEEEsetlabelwidth{$\ts{P}^{\ts{LB}}_g,\ts{P}^{\ts{UB}}_g$}\IEEEusemathlabelsep]
\item[$\ts{P}^{\ts{LB}}_g,\ts{P}^{\ts{UB}}_g$] Upper, lower bound of generator $g$
\item[$\ts{P}^{\ts{D}}_d$] Active power demand at $d$
\item[$c_g^{\ts{0}}, c_g^\ts{1}$] Start-up, Per-MW of generator $g$
\item[$c^{\ts{D,VOLL}}_\ts{d}$] Value of Lost Load cost at $d$
\item[$c_g^{\ts{offer}}, c_g^{\ts{bid}}$] Offer, bid costs of generator $g$
% \item[$P_X^{\ts{UB}}, P_X^{\ts{LB}}$] Upper, lower bound of MW constraint $X$
% \item[$U_X^{\ts{UB}}, U_X^{\ts{LB}}$] Upper, lower bound of NB constraint $X$
\item[$\ts{S}^{\ts{L}}_l$] Active power limit of line $l$
\item[$\ts{X}^{\ts{L}}_l$] Reactance of line $l$
\end{description}

\subsubsection*{Variables}
\begin{description}[\IEEEsetlabelwidth{$\ts{P}^{\ts{LB}}_g,\ts{P}^{\ts{UB}}_g$}\IEEEusemathlabelsep]
\item[$p^{\ts{G}}_\ts{g}$] Active power output of generator $g$
\item[$u_\ts{g}^\ts{G}$] Binary variable for on/off state of generator $g$
\item[$\alpha_\ts{d}^\ts{D}$] Proportion of demand met for demand $d$
\item[$p_\ts{g}^{\uparrow}, p_\ts{g}^{\downarrow}$] Positive (offer) and negative (bid) dispatch of generator $g$
\item[$\zeta^\ts{C}$] Pro-Rata \emph{Curtailment} of generators
\item[$\xi_\ts{g}^\ts{G}$] Pro-Rata \emph{Constraint} of generator $g$
\item[$\xi_\ts{k}^\ts{K}$] Pro-Rata constraint of \emph{Constraint} group $k$
\item[$\beta_{g,k}$] Binary variable modelling minimum pro-rata \emph{Constraint} of $g$ in all $k$
\item[$\delta_\ts{b}^\ts{B}$] Voltage angle at bus $b$
\item[$p^{\ts{L}}_{l}$] Active power flow in line $l$
\end{description}

%<<<<<<<<<<<<<<<<<<<
%SECTION CHANGE
%>>>>>>>>>>>>>>>>>>>

\subsection{Market Model}\label{sec:market_model}
This model represents the formation of the market merit order, therefore no network constraints are considered. Key outputs are the FPNs and generator on/off state, $p_g^{\ts{M}}$, $u_g^{\ts{M}}$ during the time period of interest.

\subsubsection{Generation}
Binary variable $u_\ts{g}^\ts{G}$ is defined for each generating unit g, to represent the generators on/off state. Generation for each unit g, $p_\ts{g}^\ts{G}$, is constrained to be between defined upper and lower bounds, $\ts{P}^{\ts{UB}}_\ts{g}$ \& $\ts{P}^{\ts{LB}}_\ts{g}$:
\begin{equation}\label{eq:market_genbounds2}
u^\ts{G}_\ts{g}\ts{P}^{\ts{LB}}_\ts{g} \le                        p^{\ts{G}}_\ts{g} \le u^\ts{G}_\ts{g}\ts{P}^{\ts{UB}}_\ts{g}
\end{equation}

\subsubsection{Power Balance}\label{subsec:market_powerbalance}
The power balance is defined over the sets of generators and demands. At this stage, the network is not considered, and the model is therefore represented as a copper-plate system.
\begin{equation}\label{eq:market_powerbalance}
    \sum_{g \in \mc{G}}p^{\ts{G}}_\ts{g} = \sum_{d \in \mc{D}}p^{\ts{D}}_\ts{d}
\end{equation}

\subsubsection{Demand}

Let $\ts{P}^{\ts{D}}_d$ denote the active power demand at demand node $d$. In some situations, the total generation may not be sufficient to meet the overall demand. In such cases, part of the demand is shed to maintain power balance. We define $\alpha_d$ as the proportion of demand served at node $d$. The demand model is then given as follows: $\alpha^\ts{D}_\ts{d} \in \left[0,1 \right], \forall d \in \mc{D}$:

\begin{equation}\label{eq:market_demand}
    p^{\ts{D}}_\ts{d} = \alpha^\ts{D}_\ts{d} \ts{P}^{\ts{D}}_\ts{d}
\end{equation}

\subsubsection{Overall formulation of the market model}
The objective function is to minimise generation and lost load costs subjects to power balance, generation and demand model. The overall formulation is given as follows:

\begin{subequations}
    \begin{align}
        \min \underbrace{\sum_{g \in \mc{G}}(c^1_\ts{g}p^{\ts{G}}_\ts{g} + c^0_\ts{g}u^\ts{G}_\ts{g})}_{\ts{Generation Cost}} + \underbrace{\sum_{d \in \mc{D}}c^{\ts{VOLL}}_\ts{d}(1-\alpha^\ts{D}_\ts{d})\ts{P}_d^D}_{\ts{Load shedding cost}}\\
        \ts{subject to~~~~~~~~~~~~~~~~~~~~~~~~~~~~~~~~~~~~~~~~} \nonumber\\
        (1-3)~~~~~~~~~~~~~~~~~~~~
    \end{align}
\end{subequations}

%<<<<<<<<<<<<<<<<<<<
%SECTION CHANGE
%>>>>>>>>>>>>>>>>>>>

\subsection{\emph{Curtailment} Model}\label{sec:curtailment_model}

This model is a surrogate for security considerations defined by the transmission system operators. The key inputs from the previous model are the market positions and on/off state for each generator, $P_g^{\ts{M}}$ and $U_\ts{g}^{\ts{M}}$ (note these are capitalised as in this model they are parameters, not variables). The key outputs are updated variables $p_\ts{g}^{\ts{S}}$ and $u_\ts{g}^{\ts{S}}$.

\subsubsection{Generation}
The constraint (\ref{eq:market_genbounds2}) remains in effect. We introduce positive variables $p^{\uparrow}_\ts{g}, p^{\downarrow}_\ts{g}$ for each unit g (offer to increase generation, bid to reduce generation), redispatch relative to $\ts{P}^{\ts{M}}_g$ is then:
\begin{equation}\label{eq:curtailment_bidoffer}
    p^\ts{G}_\ts{g} = \ts{P}^{\ts{M}}_g + p_\ts{g}^{\uparrow} - p_\ts{g}^{\downarrow}
\end{equation}

In order to model pro-rata \emph{Curtailment} of renewables, the variable $\zeta^\ts{C} \in \left[0,1\right]$ is defined to control power output 
$\forall g \in \mc{G}_{ns}$:
\begin{equation}\label{eq:curtailment_prorata}
        p^\ts{G}_\ts{g} = \ts{P}^{\ts{M}}_\ts{g} \times \zeta^{C}
\end{equation}

\subsubsection{Demand \& Power Balance}
The constraints defined in (\ref{eq:market_powerbalance}) and (\ref{eq:market_demand}) remain in place.

\subsubsection{Security constraints}
Security constraints are defined based on those in the SEMO `Weekly Operational Constraints Report' \cite{SEMO_weeklyconstraintupdate}. Only SNSP and MUON constraints are built into the model. Rate of Change of Frequency (RoCoF) and Inertia only account for circa 0.6\% of \emph{Curtailment} \cite{10253224}, and are neglected to limit model complexity.

\subsubsection{SNSP}
SNSP constraints apply to non-synchronous generators, where $\mc{G}_{ns}$ include non-synchronous generators, and $\mc{G}_{I+}$ and $\mc{G}_{I-}$ are sets of exporting and importing interconnectors respectively. This is a simplification, as pumped-hydro and BESS are not included in the model, which would also appear on the denominator as per the definition in Section \ref{sec:lit_review}. The SNSP limit is a parameter defined as $x_\ts{SNSP} \in [0,1]$. In current operations $x_\ts{SNSP} = 0.75$ \cite{SEMO_weeklyconstraintupdate}.
\begin{equation}\label{eq:SNSP}
                    \sum_{g \in \mc{G}_{ns}} p^\ts{G}_\ts{g} + \sum_{g \forall \mc{G}_{I+}} p^\ts{G}_\ts{g}\leq
        x_{SNSP} \left( \sum_{d \in \mc{D}} \ts{P}^\ts{D}_\ts{d} +
                    \sum_{g \in \mc{G}_{I-}} -p^\ts{G}_\ts{g}\right)
\end{equation}

\subsubsection{MUON}
Constraints are defined as `MW' (limit MW ouput), `MWR' (limit MW output + reserve), or `NB' (limit on/off status), and are applied to individual units or groups of units. We define $\mc{G}_\ts{X} \subset \mc{G}$ to contain each generator controlled by constraint $\ts{X}$. The constraints are defined according to the generic inequalities in (\ref{eq:MW_Examplea}).

The sets and parameters for econstraints that apply in the Irish system are defined in \ref{tab:simple_MW_constraints} and \ref{tab:simple_nb_constraints_reqs} in the Appendix. Constraints IDs reference those defined by EirGrid \& SONI in \cite{SEMO_weeklyconstraintupdate} (e.g. \texttt{S\_MWMAX\_NI\_GT} constrains output from set of generators $G_\texttt{S\_MWMAX\_NI\_GT}$ to maintain sufficient replacement reserve in Northern Ireland). Note that inter-area flow  constraints between ROI and NI (e.g \texttt{S\_MWR\_ROI}) are not included here in MUON, as they can be included in line flow limits in the DCOPF model. Some NB constraints only apply in specific circumstances based on model parameters (e.g. demand); these are defined in \ref{tab:simple_nb_constraints_reqs}. 

\begin{subequations}\label{eq:MW_Examplea}
    \begin{align}
P_\ts{X}^{\ts{LB}} \leq \sum_{g \in \mc{G}_\ts{X}}p^\ts{G}_\ts{g} &\leq P_\ts{X}^{\ts{UB}}\\
U_\ts{X}^{\ts{LB}} \leq \sum_{g \in \mc{G}_\ts{X}}u^\ts{G}_\ts{g} &\leq U_\ts{X}^{\ts{UB}}
    \end{align}
\end{subequations}

The NB (On/Off) constraints \texttt{S\_MBMIN\_CPS} and \texttt{S\_NBMIN\_MP\_NB} are enforced dependent on the variable $p^\ts{G}_\ts{g}$. To retain linearity of the model big-M formulations are used to enforce constraints. The Big-M parameters are calculated based on $\sum_{g \in \mc{G}_\ts{X}}P_\ts{g}^{\ts{UB}}$ and $\sum_{g \in \mc{G}_\ts{X}}P_\ts{g}^{\ts{LB}}$ to ensure they are as small as possible to avoid numerical solver issues.

\subsubsection*{North West Generation Constraint (\texttt{S\_MBMIN\_CPS})} 
This constraint ensures that Coolkeeragh unit C30 is online, and is enforced when the NI system demand is greater than or equal to 1,550 MW, CGT8 is unavailable and NI wind generation $<$ 450 MW. As availability is not included in the model the CGT8 requirement is neglected.
Relevant generators are defined in sub-set $\mc{G}_\ts{CPS}$. A big-M formulation is used to make the constraint conditional, introducing binary control variables $y^\ts{G}_\ts{CPS}, y_\ts{CPS}$, binary demand parameter $\ts{y}^\ts{D}_\ts{CPS} \in \left\{{0,1}\right\}$ and upper and lower big-M parameters, $\ts{M}_\ts{CPS-U}$, $\ts{M}_\ts{CPS-L}$.

\begin{subequations}\label{eq:S_MBMIN_CPS}
\allowdisplaybreaks
    \begin{align}
        \ts{y}^\ts{D}_\ts{CPS} = \begin{cases} 
              1 & \ts{if} \sum_{d \in \mc{D}_\ts{NI}} \ts{P}^\ts{D}_\ts{d} \geq 1,550 \\
              0 & \ts{otherwise}
        \end{cases}\\
        \ts{M}_\ts{CPS-U} = \sum_{g \in \mc{G}_\ts{NI,Wind}} \ts{P}^\ts{UB}_\ts{g} - 450\\
        \ts{M}_\ts{CPS-L} = 450 - \sum_{g \in \mc{G}_\ts{NI,Wind}} \ts{P}^\ts{LB}_\ts{g}\\
        450 - \sum_{g \in \mc{G}_\ts{NI,Wind}} p_g^G \leq \ts{M}_\ts{CPS-L}(y^\ts{G} _\ts{CPS})\\
        \sum_{g \in \mc{G}_\ts{NI,Wind}} p_\ts{g}^\ts{G} - 450 \leq \ts{M}_\ts{CPS-U}(1 - y^\ts{G}_\ts{CPS})\\
        y_\ts{CPS} \leq y^\ts{G}_\ts{CPS}\\
        y_\ts{CPS} \leq \ts{y}^\ts{D}_\ts{CPS}\\
        y_\ts{CPS} \geq y^\ts{G}_\ts{CPS} + \ts{y}^\ts{D}_\ts{CPS} - 1\\
        \sum_{g \in \mc{G}_\ts{CPS}} u^\ts{G}_\ts{g} \geq 1 \times y_\ts{CPS}
    \end{align}
\end{subequations}

\subsubsection*{400 kV Network Support (\texttt{S\_NBMIN\_MP\_NB})}
This constraint ensures at least one generator in $G_\texttt{S\_NBMIN\_MP\_NB}$ is operational enforced when wind generation in Republic of Ireland ($\mc{G}_\ts{ROI,W} \subset \mc{G}$) is less than 1,000 MW. Relevant generators are defined in sub-set $\mc{G}_\ts{400kV}$. A big-M formulation is used to make the constraint conditional, introducing binary variable $y^\ts{G}_\ts{400kV}$ and Big-M parameters $\ts{M}_\ts{400kV}^{\ts{UB}}$ \& $\ts{M}_\ts{400kV}^{\ts{LB}}$.

\begin{subequations}\label{eq:S_NBMIN_MP_NB}
\allowdisplaybreaks
    \begin{align}
        \ts{M}_\ts{400kV}^{\ts{UB}} = \sum_{ g \in \mc{G}_\ts{ROI,W}} \ts{P}^{\ts{UB}}_\ts{g} - 1000\\
        \ts{M}_\ts{400kV}^{\ts{LB}} = 1000 - \sum_{ g \in \mc{G}_\ts{ROI,W}} \ts{P}^{\ts{LB}}_\ts{g}\\
        1000 - \sum_{ g \in \mc{G}_\ts{ROI,W}} p_\ts{g}^\ts{G} \geq \ts{M}_\ts{400kV}^{\ts{LB}}(y^\ts{G} _\ts{400kV})\\
        \sum_{g \in \mc{G}_\ts{ROI,W}} p_\ts{g}^\ts{G} - 1000 \leq \ts{M}_\ts{400kV}^{\ts{UB}}(1 - y^\ts{G}_\ts{400kV})\\
        \sum_{g \in \mc{G}_\ts{400kV}} u^\ts{G}_\ts{g} \geq 1 \times y_\ts{CPS}
    \end{align}
\end{subequations}

\subsubsection{Overall formulation of Curtailment Model}
Bid and offer prices for the re-dispatch of generators are relative to $\ts{P}_\ts{g}^\ts{M}$ as per (\ref{eq:constraint_bidoffer}). The overall formulation is given as follows:

\begin{subequations}\label{eq:objective_curtailment}
    \begin{align}
        \min \underbrace{\sum_{g \in \mc{G}}(c_\ts{g}^\ts{offer}p_\ts{g}^{\uparrow} + c_\ts{g}^\ts{bid}p_\ts{g}^{\downarrow} + 
        c_\ts{g}^\ts{0}u_\ts{g}^\ts{G}(1-U_\ts{g}^\ts{M})}_{\ts{Redispatch Cost}}~~~~~~\nonumber\\
        + \underbrace{\sum_{d \in \mc{D}}c^{\ts{VOLL}}_d(1-\alpha^\ts{D}_\ts{d})\ts{P}_\ts{d}^\ts{D}}_{\ts{Load shedding cost}}\\
        \ts{subject to~~~~~~~~~~~~~~~~~~~~~~~} \nonumber\\
        (\ref{eq:market_genbounds2}, \ref{eq:market_powerbalance}, \ref{eq:curtailment_bidoffer} - \ref{eq:S_NBMIN_MP_NB})~~~~~~~~~~~~~~~~~
    \end{align}
\end{subequations}

%<<<<<<<<<<<<<<<<<<<
%SECTION CHANGE
%>>>>>>>>>>>>>>>>>>>

\subsection{Constraint Model}\label{sec:constraint_model}
This model represents the final dispatched position of generators on the network. The key inputs from the previous model are the re-dispatched generator set-points and on/off state for each generator, $\ts{P}_\ts{g}^{\ts{S}}$ \& $\ts{U}_\ts{g}^{\ts{S}}$. The key outputs are updated constrained generator variables, $p_\ts{g}^{\ts{C}}$ \& $u_\ts{g}^{\ts{C}}$.

\subsubsection{Generation}
The constraints defined in (\ref{eq:market_genbounds2}) remain in effect. The bid/offer constraint in (\ref{eq:curtailment_bidoffer}) is reformulated to be relative to $\ts{P}_\ts{g}^{S}$: 
\begin{equation}\label{eq:constraint_bidoffer}
    p^\ts{G}_\ts{g} = \ts{P}^{\ts{S}}_\ts{g} + p^{\uparrow}_\ts{g} - p^{\downarrow}_\ts{g}
\end{equation}

\subsubsection{Demand}
Demand constraints remain as in the previous model.

\subsubsection{Power Balance}
The power balance constraint (\ref{eq:market_powerbalance}) is replaced with (\ref{eq:constraint_powerbalance}). Let $\mc{G}_b$, $\mc{D}_b$ and $\mc{L}_b$ be sets of generators, demand and lines at bus $b$, then $\forall b \in \mc{B}$:

\begin{equation}\label{eq:constraint_powerbalance}
\sum_{g \in \mc{G}_b}p^{\ts{G}}_\ts{g} = \sum_{d \in \mc{D}_b}\ts{P}^{\ts{D}}_\ts{d}+\sum_{l \in \mc{L}_b}p^{\ts{L}}_\ts{l}
\end{equation}

\noindent
where $p^{\ts{G}}_\ts{g}$ is active power output of generator $g$, $\ts{P}^{\ts{D}}_\ts{d}$ is active power demand at $d$, and $p^{\ts{L}}_\ts{l}$ is active power flow on a line $l$.

\subsubsection{Power Flow}
AC power flow is a non-linear function of bus voltages, and computationally challenging to solve. We employ the widely used DC Optimal Power Flow (DCOPF) linearisation, that neglects reactive power and assumed a constant voltage magnitude at all buses~\cite{WoodWollenberg}. The power flow equations for both lines and transformers are defined in (\ref{eq:powerflow_KVL}), where $p^{\ts{L}}_\ts{l}$ is active power flow in line $l$, $\delta_\ts{b'}$ and $\delta_\ts{b}$ are voltage angles are the two ends of line $l$, and $\ts{X}^{\ts{L}}_\ts{l}$ is the reactance of the line $l$.

\begin{equation}\label{eq:powerflow_KVL}
p^{\ts{L}}_\ts{l} = \frac{\delta_\ts{b'} - \delta_\ts{b}}{\ts{X}^{\ts{L}}_\ts{l}}
\end{equation}

The active power flow in a line $l$ is bounded by it's capacity $\ts{S}^{\ts{L}}_{l}$
\begin{equation}\label{eq:powerflow_linelimits}
-\ts{S}^{\ts{L}}_\ts{l}\le p^{\ts{L}}_\ts{l} \le \ts{S}^{\ts{L}}_\ts{l}
\end{equation}

\subsubsection{Security constraints}
The security constraints defined in the previous section remain in place. Inter-area flows between NI and ROI are constrained by \texttt{S\_MWR\_ROI} and \texttt{S\_MWR\_NI}. These are implemented by defining $S_\ts{l}^\ts{L}$ appropriately for constraint (\ref{eq:powerflow_linelimits}).

\subsubsection{Pro-Rata \emph{Constraint} for multiple Groups} The pro-rata \emph{Curtailment} constraint (\ref{eq:curtailment_prorata}) is replaced with the following. Let $\mc{K}=\{\mc{K}_1,\mc{K}_2,\cdots,\mc{K}_n\}$ be $n$ \emph{Constraint} groups such that the $ith$ \emph{Constraint} group is defined as $\mc{K}_i \subseteq \mc{G}$. Note that the \emph{Constraint} groups may overlap and they may not be collectively exhaustive i.e. the union of all \emph{Constraint} groups may not be equal to the set of generators. Let $\xi_\ts{k}^\ts{K} \in [0,1]$ be a variable modelling the pro-rata curtailment in Constrained group $k$. We define a set $\mc{C}_g$ that contains all the \emph{Constraint} groups of each generator $g$.

$$\mc{C}_g = \{\mc{K}_i: g \in \mc{K}_i, 1\le i\le n\}$$

The pro-rata \emph{Constraint} constraint can be modelled as follows:

\begin{equation}\label{eq:constraint_minprorata}
p^{\ts{G}}_\ts{g} = \min\{\xi_\ts{k}^\ts{K}\ts{P}^{\ts{UB}}_\ts{g}:~ k \in \mc{C}_g\}
\end{equation}

(\ref{eq:constraint_minprorata}) is a discontinuous function and problematic to solve using an LP or MILP solver, so we define an equivalent formulation using binary variables is modelled as follows. Let $\beta_\ts{g,k}$ be a binary variable

\[ \beta_\ts{g,k} = \begin{cases} 
          1 & \ts{if $\xi_k$ is minimum} \\
          0 & \ts{otherwise}
       \end{cases}
\]

Noting that $\xi_\ts{k}^\ts{K} \in [0, 1]$, then $\forall g \in \mc{G}$ :
\begin{subequations}\label{eq:constraint_prorata}
\begin{align}
-\xi^\ts{G}_\ts{g} \ts{P}_\ts{g}^\ts{S} \leq   &p^{\ts{G}}_\ts{g} \le \xi^\ts{G}_\ts{g} \ts{P}^{\ts{S}}_\ts{g}\\
\xi_\ts{k}^\ts{K} - (\xi_\ts{k}^\ts{UB} - \xi_\ts{k}^\ts{LB})(1-\beta_\ts{g,k}) \leq&\xi^\ts{G}_\ts{g} \leq \xi_\ts{k}^\ts{K}\\
\sum_{k \in \mc{C}_g} \beta_\ts{g,k} &= 1
\end{align}
\end{subequations}

\subsubsection{Overall formulation of Constraint Model}
Bid and offer prices are relative to $\ts{P}_\ts{g}^{M}$ in (\ref{eq:constraint_bidoffer}), and the cost function is updated costs for turning a unit only incurred if  $\ts{U}_\ts{g}^{S} = 0$.

\begin{subequations}\label{eq:objective_constraint}
    \begin{align}
        \min \underbrace{\sum_{g \in \mc{G}}(c_\ts{g}^\ts{offer}p_\ts{g}^{\uparrow} + c_\ts{g}^\ts{bid}p_\ts{g}^{\downarrow} + 
        c_\ts{g}^\ts{0}u_\ts{g}^\ts{G}(1-U_\ts{g}^\ts{S})}_{\ts{Redispatch Cost}}~~~~~~\nonumber\\
        + \underbrace{\sum_{d \in \mc{D}}c^{\ts{VOLL}}_\ts{d}(1-\alpha^\ts{D}_\ts{d})\ts{P}_d^D}_{\ts{Load shedding cost}}\\
        \ts{subject to~~~~~~~~~~~~~~~~~~~~~~~} \nonumber\\
        (\ref{eq:market_genbounds2},
        \ref{eq:market_demand},
        \ref{eq:SNSP} - \ref{eq:S_NBMIN_MP_NB}, \ref{eq:constraint_bidoffer} - \ref{eq:constraint_prorata})~~~~~~~~~~~~~~~~~
    \end{align}
\end{subequations}

%SECTION CHANGE
%>>>>>>>>>>>>>>>>>>>

\subsection{Allocation of \emph{Curtailment} \& \emph{Constraint} Volumes}\label{sec:allocation_of_volumes}
Volumes of Surplus, \emph{Curtailment} and \emph{Constraint} are determined post modelling. This aligns with the TSO methodology in \cite{EirGrid_flagging_methodology}, where dispatch actions are flagged by the control room systems \emph{ex-ante}.

Surplus is the difference between the maximum potential output of non-synchronous generation and their market position:
\begin{equation}\label{eq:vol_surplus}
    \ts{V}^{\ts{Surplus}} = \sum_{g \in \mc{G}_{ns}} \left(P_\ts{g}^{\ts{UB}} - \ts{P}_\ts{g}^{\ts{M}}\right)
\end{equation}

If the SNSP constraint is active, the as suggested in \cite{EirGrid_flagging_methodology} is is assumed that all curtailment to meet the limit is flagged as SNSP. $\ts{V}^{\ts{SNSP}}$ is therefore defined as:
\begin{subequations}\label{eq:vol_SNSP}
    \begin{align}
        \Delta^\ts{SNSP} =\underbrace{\sum_{g \in \mc{G}_{ns}} \ts{P}^\ts{M}_\ts{g} +                        \sum_{g \in \mc{G}_{I+}} \ts{P}^\ts{M}_\ts{g}}_{\ts{Non-Synchronous Supply}}~~~~~~~~~~~~~~~~~~~~~~~~~~~ \nonumber\\
        - \underbrace{x_\ts{SNSP} \left( \sum_{d \in \mc{D}} \ts{P}^\ts{D}_\ts{d} +
    \sum_{g \in \mc{G}_{I-}} -\ts{P}^\ts{M}_\ts{g}\right)}_{\ts{Net Loads}}\\
        \ts{V}^{\ts{SNSP}} = \max(0, \Delta^\ts{SNSP})
    \end{align}
\end{subequations}

The volume of MUON \emph{Curtailment} is defined as reductions in non-synchronous generation to secure the network beyond $\ts{V}^{\ts{SNSP}}$. 
\begin{equation}\label{eq:vol_MUON}
    \ts{V}^{\ts{MUON}} = \sum_{g\in\mc{G}_{ns}}\left(\ts{P}^{\ts{M}}_\ts{g} - \ts{P}^{\ts{S}}_\ts{g}\right) - \ts{V}^{\ts{SNSP}}
\end{equation}

% \subsection{\emph{Constraint}}
Network constraints are introduced in the final UC-DCOPF model. The amount of \emph{Constraint} can then be calculated as:
\begin{equation}\label{eq:vol_constraint}
    \ts{V}^{\ts{TC}} = \sum_{g\in\mc{G}_{ns}}\left(\ts{P}^{\ts{S}}_\ts{g} - p^{\ts{C}}_\ts{g} \right)
\end{equation}

%<<<<<<<<<<<<<<<<<<<
%SECTION CHANGE
%>>>>>>>>>>>>>>>>>>>

\section{Results}\label{sec:test}
The model described in the previous sections is demonstrated on a reduced illustrative example, and then on a representative Irish model. The formulation is coded in python using pyomo, with Gurobi as the default solver \cite{oats_curtailment}.

\subsection{An Illustrative Test Case}
\begin{figure}[h]
    \centering
    \includegraphics[width=1\linewidth]{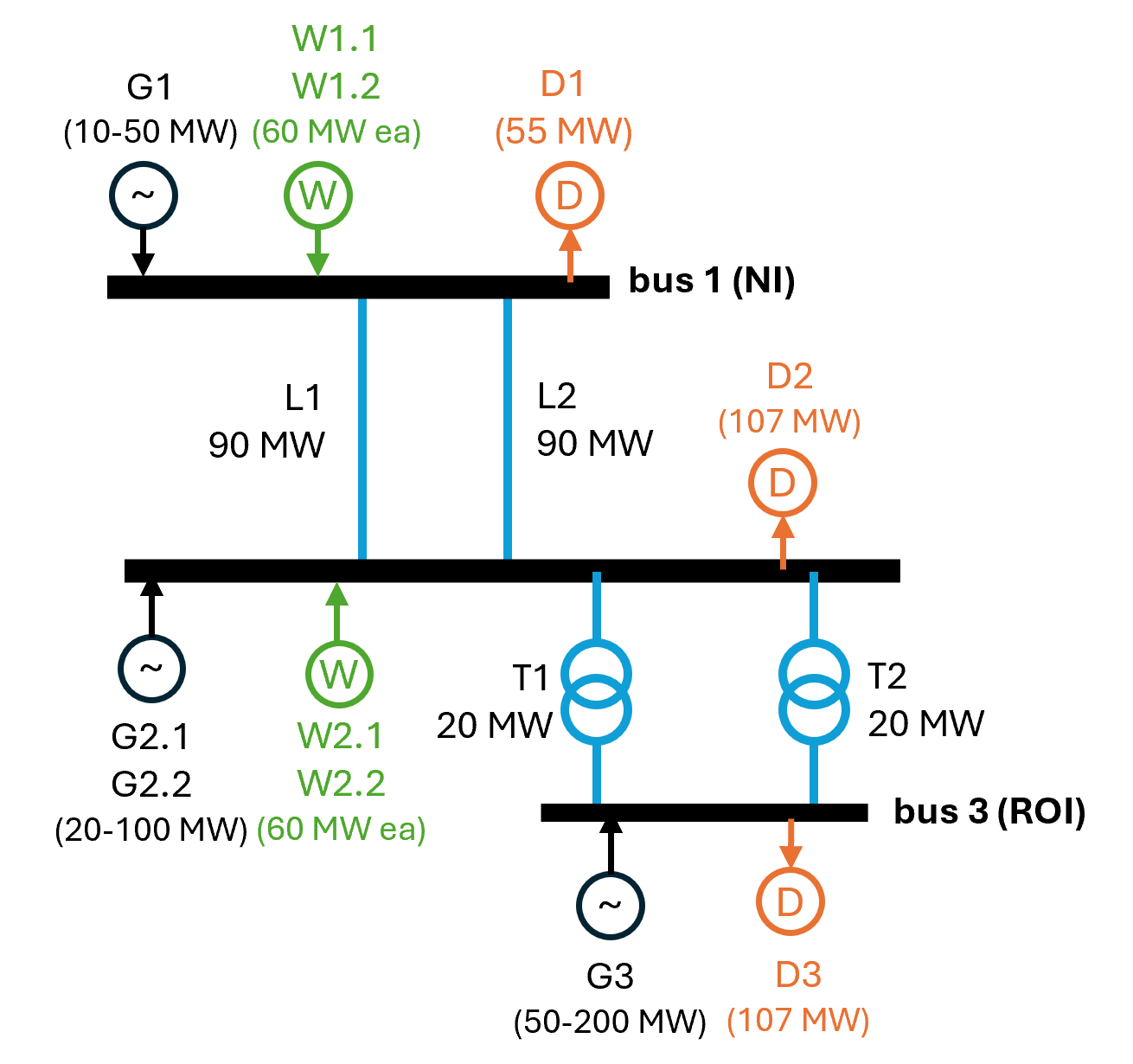}
    \caption{Illustrative Test Case (Original)}
    \label{fig:illustrative_testcase}
\end{figure}

To test model behaviour a snapshot test-case was used, based on one developed by RES, who funded development of an earlier version of a `\emph{Constraint} tool'. The diagram in see Fig. \ref{fig:illustrative_testcase} shows generator availability, demand and line/transformer ratings. $c_\ts{g}^0$ \& $c_\ts{g}^1$ were set higher for synchronous generators than wind, and $c_\ts{g}^\ts{bid}$ set lower. The MUON constraints shown in the Table \ref{tab:illustrative_example_MUON_cosntraints} were active. \emph{Constraint} groups are defined as $\mc{K}_1 = \left\{W1.1, W1.2\right\}, \mc{K}_2 = \left\{W1.2, W2.1\right\}, \mc{K}_3 = \left\{W2.1, W2.2\right\}$. Note that to demonstrate pro-rata curtailment in overlapping constraint groups $W1.1 \in \left\{\mc{K}_1,\mc{K}_2\right\}$.
\begin{table}[h]
    \centering
    \caption{Illustrative Example MUON Constraints}
    \label{tab:illustrative_example_MUON_cosntraints}
    \begin{tabular}{c|c|c}
        \toprule
        ID & $\mc{G_X}$ & Notes\\
        \midrule
        \texttt{S\_MWMAX\_NI\_GT}  &  $\left\{G1\right\}$   & $\sum_{g\in\mc{G}_x}{p^G_g} \leq 80$ \\
        \texttt{S\_NBMIN\_DUB\_L2} &  $\left\{G2.1, G2.2\right\}$   & $\sum_{g\in\mc{G}_x}{u_g} \geq 1$ \\
        \bottomrule
    \end{tabular}
\end{table}

\begin{figure}
    \centering
    \includegraphics[width=1\linewidth]{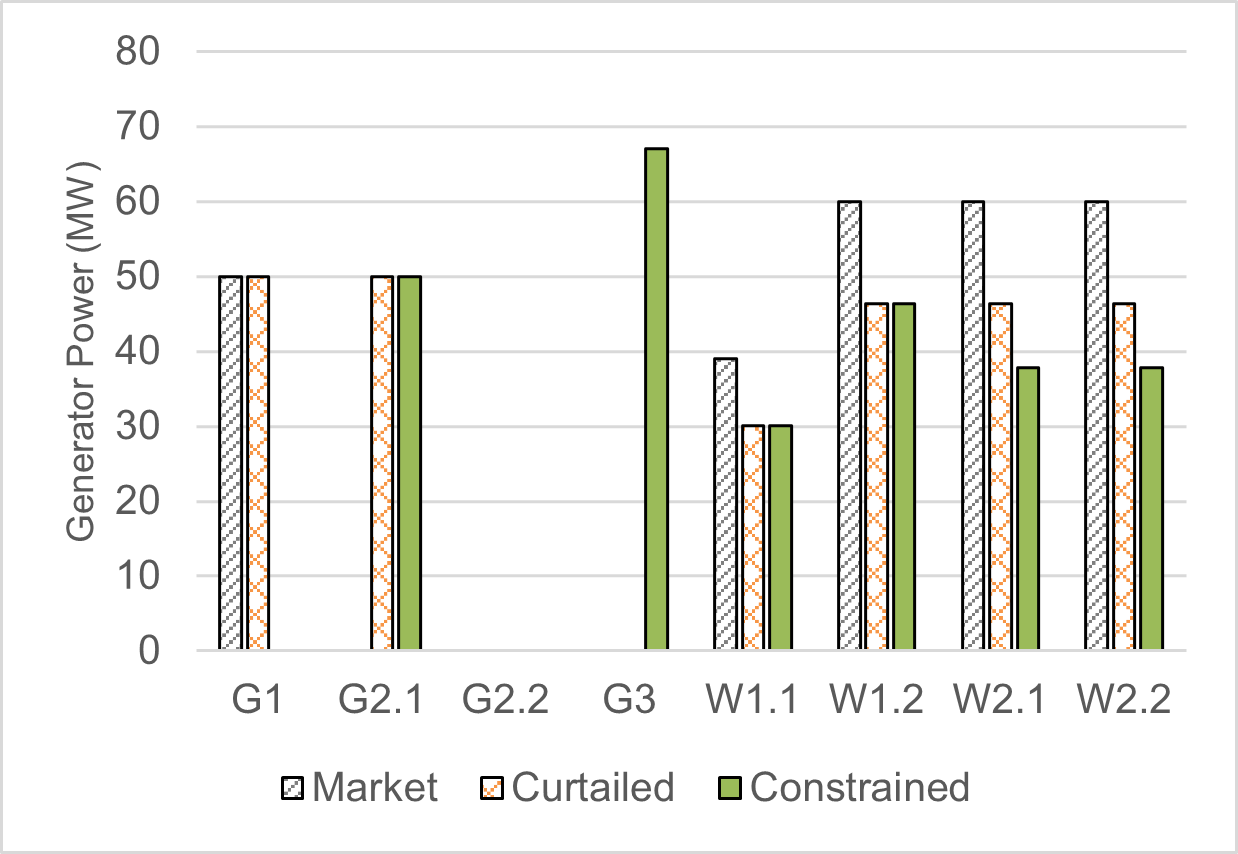}
    \caption{Illustrative Test Results}
    \label{fig:illustrative_test_result}
\end{figure}

The results (see Fig.~\ref{fig:illustrative_test_result}) in MW$\cdot\Delta t$ ($\Delta t$ being the timestep length, which is undefined for this test) are:
\begin{itemize}
    \item $V^\ts{Surplus} = 21.00$
    \item $V^\ts{SNSP} = 17.25$
    \item $V^\ts{MUON} = 32.75$
    \item $V^\ts{TC} = 17.00$
\end{itemize}

The market allocates wind where possible, however due to the minimum generation requirement of G1, wind from W1.1 is displaced, leading to surplus. The inclusion of minimum generation levels in the market model reflects the ability to submit `Minimum Acceptance Volumes' for Price-Quanity pairs in the market \cite{SEMO_MarketProcedures}. The impact of this behaviour will be small on a realistic Ireland network model, as the minimum generation levels of synchronous generators will be smaller relative to the total system demand, and non-synchronous generation unlikely (at present) to approach such high proportions of available generation.

Wind is then \emph{curtailed} pro-rata across the system due to SNSP \emph{Curtailment}, and MUON \emph{Curtailment}. In the \emph{Constraint} model wind generation in $\mc{K}_{1}$ at bus 2 is dispatched down 'pro-rata', as demand at D3 is 107 MW, however only 40 MW of transmission capacity exists. Despite W1.2 belonging to a shared pro-rata group with W2.1, it is not dispatched down to the same level, demonstrating that the formulation in \ref{eq:constraint_prorata} achieves the desired outcome.

\subsection{All-Island Testcase}
A representative transmission system network was built for the purposes of testing the tool, including busses (446), lines (586), transformers (183) and generators (288). The construction of the model requires significant manual review of the data to address inaccuracies and inconsistencies, and is in a developed draft stage \cite{Ireland_network_model}, but is sufficient to test the operation of the model. The main source of the data for this is the All-Island Ten Year Transmission Forecast Statement 2022 (TYTFS) \cite{AllIslandTYTFS2022}, with further data taken from the SEMO Market Data Portal \cite{SEMO_MarketData} and EirGrid Quarterly System and Renewable Data Reports \cite{EirGrid_Reports}. Interconnector flows are taken from the GB Balancing Mechanism Reporting Service (BMRS) \cite{BMRS}.

\begin{figure}[h]
    \centering
    \includegraphics[width=1\linewidth]{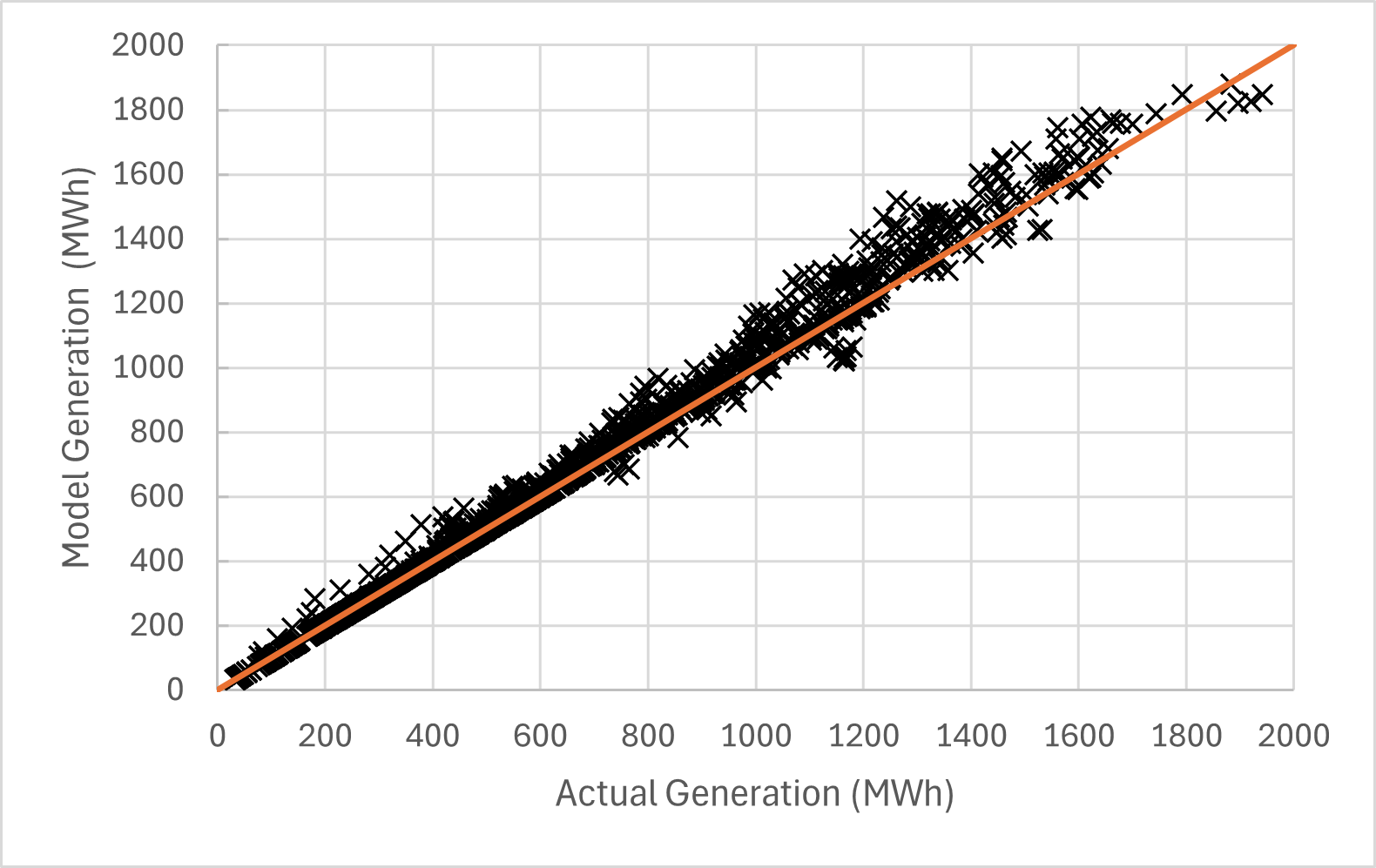}
    \caption{All-Island Testcase (July 2025) - Wind \& Solar Generation, Modelled vs. Actual}
    \label{fig:results_09_genvgen}
\end{figure}

\begin{figure}[h]
    \centering
    \includegraphics[width=1\linewidth]{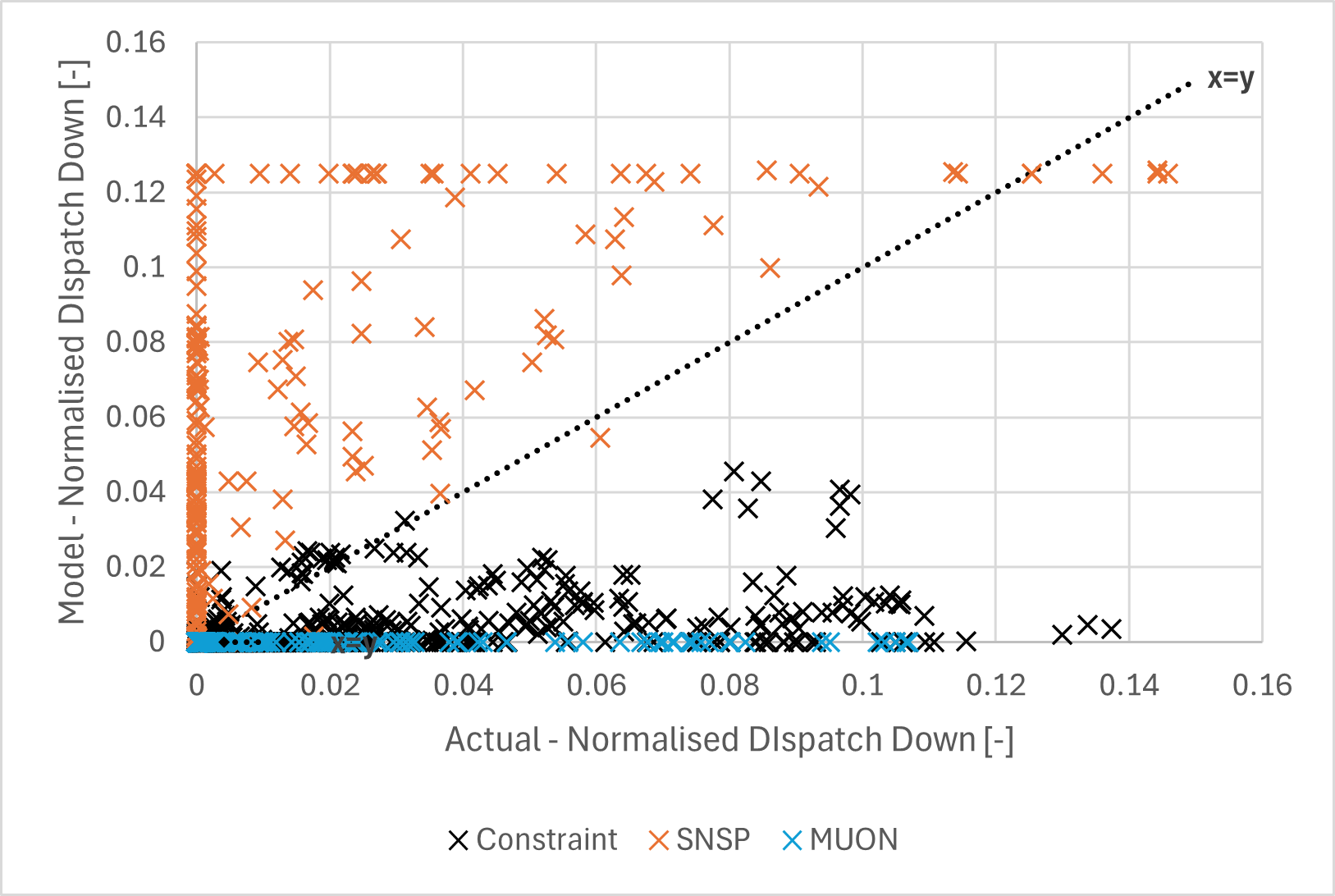}
    \caption{All-Island Testcase (July 2025) - Dispatch Down, Modelled vs. Actual}
    \label{fig:results_09_timewise}
\end{figure}

A test run for the month of July 2025 at half-hourly resolution (1,488 periods) was carried out on a Dell XPS 15 laptop with an Intel(R) Core(TM) i7-10750H CPU @ 2.60GHz and 32 GB of RAM solving in 32 minutes. The results are compared to actual dispatch down data from EirGrid \cite{EirGrid_Reports}. It can be seen in Fig.~\ref{fig:results_09_genvgen} that the actual output of wind \& solar is a reasonable estimate, but generally overestimated. It is clear however from Fig.~\ref{fig:results_09_timewise} that the model under estimates total dispatch down, \emph{Constraints} and MUON, while overestimating SNSP \emph{Curtailment}. The under prediction of \emph{Constraints} is an expected limitation of this model, as N-1 security constraints are not considered and the network model still requires validation. 

Calculating the real world actual pre \emph{curtailment} SNSP levels and volumes of SNSP \emph{Curtailment} from  \cite{EirGrid_Reports}, the SNSP volumes are lower than expected using the logic in \ref{eq:vol_SNSP}. This could be caused by demand increases from pumped-storage or BESS charging, differences in the  tagging/flagging of dispatch actions relative to the model logic (e.g. allocation to MUON or \emph{Constraints}), or the lack of time-coupling to represent minimum up-time of synchronous generators and dispatch actions from previous periods that reduce the set-point below where the SNSP constraint is bounding.  This appears to be a key driver in over-estimates of SNSP \emph{Curtailment}, such that MUON volumes are never the limiting constraint. Implementing more detailed flagging logic, and implementing time-coupling are key areas of future work.

%<<<<<<<<<<<<<<<<<<<
%SECTION CHANGE
%>>>>>>>>>>>>>>>>>>>

\section{Discussion \& Conclusions}\label{sec:conclusion}
The Irish electricity system involves complex market arrangements and Transmission System Operator (TSO) operations. While much information is available about how these systems function, there is limited technical detail on the exact commitment and dispatch formulations used by the TSO.

This paper presents a replicable mathematical modelling formulation that captures Ireland-specific security constraints, including a novel pro-rata dispatch method for generators within multiple constraint groups. A representative Irish transmission network model was developed~\cite{Ireland_network_model}, and operational data for July 2025 were assembled to perform model validation. The validation results provide confidence in the accuracy of the proposed approach, while also revealing its limitations.

The use of the model is dependent on sufficiently accurate input data, such as generator parameters (operating ranges, costs), interconnector flows, and network physical and topological characteristics. Cost challenges can be overcome by assuming cost values relative to type (e.g $c^{1}_\ts{gas}<c^{1}_\ts{coal}$) to estimate market formation. An area of ongoing work is the further development of a model of the Irish transmission network.

The authors acknowledge that the simple model is limited by the lack of time-coupling and (N-1) contingencies, leading to deviations from real world results when tested on the All-Island Testcase. Introducing time-coupling would ensure detailed generator operation and market behaviour is captured, such as ramping rates, dwell points, and minimum up/down times. The capability to time-couple dispatch instructions would also allow for more realistic modelling of the system.

Incorporating (N-1) contingency analysis would provide more realistic transmission system availability and constraint volumes; and having an understanding of interconnected regions would ensure market/price formation is more realistic. This reformulation would result in a Security Constrained Unit Commitment (SCUC) Optimal Power Flow, increasing computation time and complexity relative to the model presented here, while introducing further data requirements for accurate generator parameters.

Interconnectors are currently modelled as generators with a $P_g^{LB} < 0$ with static marginal costs, while in reality the volume and direction of flow will depend on the relative cost of the coupled markets (GB and EU). The current model does not model storage (e.g. pump hydro, or grid scale battery), which are increasingly providing flexibility services, and where paired with generators could reduce curtailment by shifting export times. 

In summary, this paper has articulated the requirements of a model of the Irish test system and formulated a simplified representation, providing an important starting point for further open source model development that can support analysis of \emph{Curtailment} and \emph{Constraint} of renewables in Ireland.

\section*{Acknowledgements}\label{sec:ack}
The research leading to these results has received funding from Renewable Energy Systems Limited (RES) through the Innovation Project under Grant agreement number 20125. The authors would also like to thank Manuel Hurtado of EirGrid for providing time and advice.

\bibliographystyle{IEEEtran}
\bibliography{biblio}{}

\appendices
\section{Security Constraints in Irish Electricity System}\label{appen:SecurityConsIreland}

Constraints IDs reference those defined by EirGrid \& SONI in \cite{SEMO_weeklyconstraintupdate} (e.g. \texttt{S\_MWMAX\_NI\_GT} constrains output from set of generators $G_\texttt{S\_MWMAX\_NI\_GT}$ to maintain sufficient replacement reserve in Northern Ireland)

\begin{table}[h]
    \centering
    \caption{Security MW Constraints}
    \label{tab:simple_MW_constraints}
    \begin{tabular}{c|c|c|c}
        \toprule
        ID & $\mc{G_X}$ & $\ts{P}_\ts{g}^{\ts{LB}}$ & $\ts{P}_\ts{g}^{\ts{UB}}$\\
        \midrule
        \texttt{S\_MWMAX\_NI\_GT} & $\mc{G}_{NI\_GT}$   & - & 272\\
        \texttt{S\_MWMIN\_EWIC}   & $g_{EWIC}$   & -526 & -  \\
        \texttt{S\_MWMAX\_EWIC}   & $g_{EWIC}$ &      &  504\\
        \texttt{S\_MWMIN\_MOYLE}  & $g_{MOYLE}$  &  -410 & \\
        \texttt{S\_MWMAX\_MOYLE}  &  $g_{MOYLE}$  &    &   441\\
        \texttt{S\_REP\_ROI}      & $\mc{G}_{ROIGT}$ & - & [1] \\
        \texttt{S\_MWMAX\_ROI\_GT} &  &  &\\
        \texttt{S\_MWMIN\_CRK\_MW} & $\mc{G}_{CRK\_MW}$ & 0 & \\
        \texttt{S\_MWMAX\_CRK\_MW} &  &  & 1,370\\
        \texttt{S\_MWMIN\_STH\_MW} & $\mc{G}_{STH\_MW}$ & 0 & \\
        \texttt{S\_MWMAX\_STH\_MW} &  &   &1,835\\ 
        \midrule
        \multicolumn{4}{l}{[1] MWR Limit: $\sum_{g \in \mc{G}_{ROIGT}} \ts{P}^{\ts{UB}}_g - 325$}\\
        \bottomrule
    \end{tabular}
    \label{tab:placeholder}
\end{table}

\begin{table}[h]
    \centering
    \caption{Security NB (On/Off) Constraints}
    \label{tab:simple_nb_constraints}
    \begin{tabular}{c|c|c|c}
        \toprule
        ID & $\mc{G_X}$ & $\ts{U}_\ts{g}^{LB}$ & $\ts{U}_\ts{g}^{UB}$\\
        \midrule
        \texttt{S\_NBMIN\_MINNIU} & $\mc{G}_{MINNIU}$   & 3 & -\\
        \texttt{S\_NBMIN\_MINNI3} & $\mc{G}_{MINNI3}$   & 1 & -\\
        \texttt{S\_NBMIN\_ROImin} & $\mc{G}_{ROImin}$   & 4 & -\\
        \texttt{S\_NBMIN\_DubNB}  & $\mc{G}_{DubNB}$    & 1 & -\\
        \texttt{S\_NBMIN\_DubNB2} & $\mc{G}_{Dub\_NB2}$ & 2 & -\\
        \texttt{S\_NBMIN\_DUB\_L1}& $\mc{G}_{DUB\_L1}$  & 3 & -\\
        \texttt{S\_NBMIN\_DUB\_L2}& $\mc{G}_{DUB\_L2}$  & 1 & -\\
        \texttt{MP5\_NB}          & $\mc{G}_{MPX}$      & - & 1\\
        \bottomrule
    \end{tabular}
\end{table}

\begin{table}[h]
    \centering
    \caption{Security NB (On/Off) Constraints Logic}
    \label{tab:simple_nb_constraints_reqs}
    \begin{tabular}{c|c|c|c}
        \toprule
        ID & Parameter & Operator & Limit\\
        \midrule
        \texttt{S\_NBMIN\_DUB\_L1}& $\sum_{d \in \mc{D}_{ROI}} P^{\ts{D}}_d$  & $\geq$ & 4,000\\
        \texttt{S\_NBMIN\_DUB\_L2}& $\sum_{d \in \mc{D}_{ROI}} P^{\ts{D}}_d$  & $\geq$ & 4,700\\
        \bottomrule
    \end{tabular}
\end{table}

\end{document}